\documentclass[review,hidelinks,onefignum,onetabnum]{siamart250211}

\nolinenumbers

\usepackage{lipsum}
\usepackage{amsfonts}
\usepackage{graphicx}
\usepackage{epstopdf}
\usepackage{amssymb,amsfonts,amsmath,latexsym,dsfont}
\usepackage{mathtools}
\usepackage{mathrsfs}
\usepackage[algo2e,ruled,vlined]{algorithm2e}
\usepackage{setspace}
\usepackage{subfigure}
\usepackage{fancyhdr}
\usepackage{wrapfig}

\usepackage{booktabs}

\usepackage[normalem]{ulem}

\usepackage{boxedminipage}
\usepackage{pgf,pgfarrows} 
\usepackage{subfigure} 

\newcommand{\FrameboxA}[2][]{#2}
\newcommand{\Framebox}[1][]{\FrameboxA}

\newcommand{\mc}[3]{\multicolumn{#1}{#2}{#3}}

\newcommand{\im}{\textit{\i}}

\newcommand{\bfe}{{\bf e}}
\newcommand{\bfb}{{\bf b}}

\newcommand{\bfx}{ {\bf x}}

\newcommand{\bfu}{{\bf u}}

\newcommand{\bfr}{{\bf r}}

\renewcommand{\theequation}{\arabic{section}.\arabic{equation}}

\definecolor{darkblue}{rgb}{0.08, 0.15, 0.48}

\newtheorem{remark}{Remark}[section]

\graphicspath{{images/}{./}}

\ifpdf
  \DeclareGraphicsExtensions{.eps,.pdf,.png,.jpg}
\else
  \DeclareGraphicsExtensions{.eps}
\fi

\headers{Real-Shifted Coarse Grid Correction for Helmholtz}{R. Yovel and E. Treister}

\title{Scalable Multigrid Solver for the Helmholtz Equation: Real-Shifted Coarse Grid Correction  \thanks{Corresponding author: Rachel Yovel. \funding{This research was supported by The Israel Science Foundation (grant No. 656/23). RY is supported by the Ariane de Rothschild scholarship and by Kreitman High-tech scholarship. The authors also thank the Lynn and William Frankel Center for Computer Science at BGU.}}}

\author{Rachel Yovel\thanks{Institute for Interdisciplinary Computational Science, the Stein Faculty of Computer and Information Science, Ben-Gurion University of the Negev, Beer-Sheva, Israel.
  \email{yovelr@bgu.ac.il, erant@bgu.ac.il}}
\and Eran Treister$^\dag$}

\usepackage{amsopn}

\begin{document}

\maketitle

\begin{abstract}
We present a convergent and scalable multigrid solver for high-frequency Helmholtz equations.
Standard multigrid methods do not converge for high-frequency Helmholtz problems, and a common cure is adding a complex shift and using the shifted operator as a preconditioner.
Nevertheless, the complex shift prevents scalability.
In this work we present a new method that achieves scalable convergence of a 3-level cycle without a complex shift.
Our key idea is real-shifting the coarsest grid Galerkin operator, to correct the numerical dispersion between the grids.
We show that this real-shifted coarse grid correction leads to a scalable 3-level method, for  problems with 12 grid points per wavelength on the fine grid, and a convergent cycle with very few iterations for 11 grid points per wavelength, using standard point-smoothers.
For problems with 10 grid points per wavelength, our method combined with a modest complex shift outperforms the standard complex shifted Laplacian method by an order of magnitude. 
We demonstrate wavenumber independent convergence for  heterogeneous geophysical media in 2D and 3D.
\end{abstract}



\section{Introduction}\label{sec:intro}


The Helmholtz equation models wave propagation in the frequency domain. 
In non-dissipative acoustic media, it is given by
\begin{equation} \label{eq:acousitcHelm}
-\Delta p(\vec{x}) - k^2(\vec{x})p(\vec{x}) = q(\vec{x}) \qquad \vec{x}\in\Omega
\end{equation}
where $p$ is the Fourier transform of the wave's pressure field, $k(\vec x) = \omega \kappa(\vec x)$ is the wavenumber, comprised of the angular frequency $\omega$ and the slowness $\kappa$ (the inverse of the wave velocity). Finally, $q$ is the source term.
Eq. \eqref{eq:acousitcHelm} is typically discretized using finite-differences on a regular grid, with absorbing boundary conditions (ABC) \cite{engquist1977absorbing} or perfectly matched layers (PML)
\cite{berenger1994perfectly,   harari2000analytical, rabinovich2010comparison} to mimic propagation in an open domain. 
For high-frequency Helmholtz problems, the resulting linear system is indefinite and hence solving the resulting problem is very challenging numerically.
High-frequency waves require fine meshes (typically, 10 grid points per wavelength are taken in many solution methods \cite{bayliss1985accuracy, haber2011fast}), making the linear system very large.
A sought property in many solution methods is \emph{wavenumber independent convergence}. When keeping the number of grid points per wavelength constant, this property translates to having a nearly constant iteration count, regardless of the grid size.
Suggested preconditioners include mainly domain decomposition \cite{benamou1997domain, gander2013domain, stolk2013rapidly, chen2016robust} and multigrid methods, such as the complex-shifted Laplacian method \cite{erlangga2006novel},
smoothed aggregation algebraic multigrid methods \cite{olson2010smoothed, tsuji2015augmented},
a multigrid solver to the Helmholtz equation with a point source based on travel time and amplitude \cite{treister2019multigrid},
 wave-ray methods \cite{livshits2014scalable}, 
and algebraic multigrid solver with a new coarse correction \cite{falgout2025toward}.
Yet, existing methods of more than two levels are not known to achieve wavenumber independent convergence, unless using expensive wave-ray formulations, as in \cite{livshits2014scalable}.
Recently, deep-learning based multigrid methods were also suggested for the Helmholtz equation \cite{azulay2022multigrid, cui2025neural},
 however, they require expensive training.

The main challenge in solving the Helmholtz equation using multigrid methods lies in the discrepency between the fine grid operator and its coarse approximations in a spectral sense. 
In \cite{elman2001multigrid}, it is shown that the error introduced by the coarse grid correction in a two-grid method is proportional to
\begin{equation}\label{eq:1-lambdaRatio}
1-\frac{\lambda_h}{\lambda_H}
\end{equation}
where $\lambda_h$ is a fine eigenvalue and $\lambda_H$ is the eigenvalue corresponding to the coarsening of the corresponding eigenvector.
Hence, when an eigenvalue undergoes a sign change during the coarsening step, the multigrid cycle diverges.
The remedy suggested in \cite{elman2001multigrid} is adding a complex shift, which, if sufficiently large, promises convergence.
Based on this idea, the well-known complex-shifted Laplacian multigrid preconditioner (CSLP) \cite{erlangga2006novel} suggests using a shifted version as a preconditioner for the original system, solved by multigrid. 
Nevertheless, CSLP is bound to lose scalability: 
it requires a shift of $O(k^2)$ for convergence \cite{cocquet2017large},
while more than $O(k)$  prevents scalability \cite{gander2015applying}.

\emph{Numerical dispersion} is the misalignment between the continuous and discrete wave speeds, caused by the discretization scheme. 
In the context of multigrid, 
discrepancies between $\lambda_h$ and $\lambda_H$ in \eqref{eq:1-lambdaRatio} manifest as numerical dispersion (phase shifts), which
can accumulate to yield opposite phases over multiple wavelengths and degrade convergence.
Even a small phase shift, after accumulation can yield opposite phases.
This phenomenon is referred to as the \emph{pollution} effect.
It was observed \cite{stolk2014multigrid, treister2019multigrid, cocquet2021closed, yovel2024lfa, stolk2025two} that lower dispersion error improves multigrid convergence.
While high-order discretizations for the Helmholtz equation improve multigrid convergence, their ability to correct dispersion is limited and complex shifts are still needed \cite{umetani2009multigrid, yovel2025vanka}.

In \cite{jo1996optimal}, 
a family of compact dispersion-correcting discretizations was suggested. 
These methods require re-optimization of 9 parameters for each configuration, and remain only second order accurate.
The recent assymptotic method \cite{cocquet2024asymptotic} suggests optimizing only a real shift of the wavenumber while keeping the original order of accuracy, though, the dispersion correction is not optimal for extremely small numbers of points per wavelength.
While in 1D the numerical dispersion can be completely eliminated by a real shift \cite{ernst2013multigrid}, 
in more than 1D, it is commonly agreed that numerical dispersion cannot be eliminated and the pollution effect is bound to appear after some wavelengths. 
As a result, reducing numerical dispersion in higher dimensions primarily serves to delay the onset of the pollution effect rather than eliminate it, which can still be beneficial for enabling multigrid convergence and reducing the need for complex shifts.

Some works suggest improved multigrid cycles using dispersion corrected re-discretization methods.
In particular, \cite{stolk2014multigrid} propose a two-level method in which the weights of each stencil are optimized to minimize numerical dispersion. 
A small complex shift was also applied to ensure convergence. 
While this approach demonstrated improved convergence starting from approximately 3.5 grid points per wavelength on the coarsest grid, it was restricted to two-level cycles.
Similarly, \cite{cocquet2021closed} studied dispersion-corrected  discretizations, demonstrating improvements in multigrid convergence. 
Their experiments, conducted without any complex shift, showed convergence with 3 grid points per wavelength on the coarsest level. 
However, achieving convergence required 10--20 expensive Kaczmarz relaxations per level, and the method was tested on relatively small 2D grids only. 
By presenting a Galerkin-coarsened dispersion correction method, we show that only modest smoothing --- one pre- and one post relaxation of damped Jacobi --- is sufficient, even for larger 2D and 3D grids.

In this work, we introduce a real-shifted coarse grid correction (RS-CGC) for the Helmholtz equation.
Our method achieves superior performance with no complex shift, while saving the need to re-discretize and optimize finite difference weights:
we interweave a real shift into the Galerkin coarsening method, combined with high-order disctetization and intergrid operators.
Since the three-level cycle that we suggest is applied to the original operator, and not to a shifted version, our method can also serve as a solver, and not only as a preconditioner for a Krylov method.

To tune the real shift and validate the results, we develop a specialized grid-to-grid dispersion analysis, focusing on the misalignment of the first and the third grid.
Unlike local Fourier analysis (a predictive tool for convergence of multigrid methods), which couples smoothing and coarse grid correction, our analysis assesses the optimal shifts and limitations of the coarse grid correction regardless of the smoother.
Unlike standard dispersion analysis that compares a discrete operator to a continuous one, our method is designed to address multigrid performance.

The resulting real shifts save the need of a complex shift and enable a scalable multigrid method, which can serve as a preconditioner inside a Krylov method or act as a scalable solver, for as few as 12 grid points per wavelength (3 points on the coarsest grid).
The method is still rapidly convergent for 11 grid points per wavelength, with only 2.75 grid points on the coarsest level.
For 10 grid points per wavelength, we suggest a combination of RS-CGC with CSLP which gains a speedup of an order of magnitude over classical CSLP.
Even though our analysis does not cover spatially dependent wavenumbers, the use of the same shifts that correspond to the constant wavenumber case is efficient for the spatially dependent case as well. 
In particular, we show similar scalability properties for  heterogeneous geophysical media in 2D and 3D.

The paper is organized as follows: 
In Section \ref{sec:mg} we briefly give multigrid preliminaries.
In Section \ref{sec:main} we describe our method in detail and in Section \ref{sec:analysis} we derive our grid-to-grid analysis, choose the real shifts and show numerical examples for the theory.
In Section \ref{sec:numerical} we demonstrate the achievements by numerical experiments, and finally in Section \ref{sec:conclusion} we provide concluding remarks and notes about future work.

\section{Multigrid preliminaries}\label{sec:mg}

Multigrid methods \cite{brandt1977multi} form a family of iterative methods for solving linear systems obtained from PDE discretizations, particularly systems of the form $H_h\bfx=\bfb$.
For elliptic equations, multigrid methods are very efficient.
Their efficiency stems from the interaction of two complementary processes. 
\emph{Smoothing} --- the application of relaxation methods, such as damped Jacobi or Gauss-Seidel --- is used to damp high-frequency error components.
To complement it, the \emph{coarse-grid correction} --- using information from an exact solution on a coarser grid to correct the error ---
 is employed to remove the remaining low-frequency error modes.
Two operators are introduced to translate the information between the grids:
the \emph{restriction} $I_h^{2h}$, that translates the residual from fine to coarse grid, and the \emph{prolongation}, which interpolates the corrected error back to the fine grid.
We summarize this two-grid process in Algorithm \ref{alg:TwoCycle}. 
Recursive application of Algorithm \ref{alg:TwoCycle} results in the multigrid cycle. 
Using one recursive call, we obtain the multigrid V-cycle, and two recursive calls gives the multigrid W-cycle.
For more details, see \cite{briggs2000multigrid, trottenberg2000multigrid} and references therein.

\begin{algorithm}
\DontPrintSemicolon
\KwSty{Algorithm: $\bfu\leftarrow TwoGrid(H_h,\bfb,\bfx).$}\;
\begin{enumerate}\Indm
\item Apply pre-relaxations: $\bfx \leftarrow Relax(H_h,\bfx,\bfb)$\;
\item Compute and restrict the residual $\bfr_c = I_h^{2h}(\bfb - H_h\bfx)$.
\item Compute $\bfe_{2h}$ as the solution of the coarse-grid problem $H_{2h}\bfe_{2h}=\bfr_{2h}$.
\item Apply coarse grid correction: $\bfx \leftarrow \bfx + I_{2h}^h\bfe_{2h}$.
\item Apply post-relaxations: $\bfx \leftarrow Relax(H_h,\bfx,\bfb)$.
\end{enumerate}
\caption{Two-grid cycle.}
\label{alg:TwoCycle}
\end{algorithm}

For high-frequency Helmholtz problems, standard multigrid methods  fail to converge. 
Standard relaxation methods (such as damped Jacobi) are not stable smoothers as they diverge when the number of relaxation cycles increases, and the coarse grid correction fails because of the sign switching in \eqref{eq:1-lambdaRatio} after a mode undergoes coarsening. 
An added complex shift improves not only the coarse grid correction, but also stabilizes the smoother as the operator becomes closer to diagonal dominance.
For Poisson-like problems, V-cycles and W-cycles performs similarly in terms of iteration count (and hence V-cycles are cheaper in computational time), since the coarse grid correction is nearly ideal.
However, for Helmholtz problems, it is very common to use W-cycles since the coarse grid correction (even after adding a shift) struggles to converge.
We present a convergent multigrid cycle for high-frequency Helmholtz problems, which converges for both V- and W-cycle, despite the use of unstable point smoothers and very few smoothing cycles.
We leave the improvement of the smoother for future work.

In matrix form, the two-grid operator described in Algorithm \ref{alg:TwoCycle} can be written as
\begin{equation}\label{eq:TG}
TG = S^{\nu_2}\left(I - I_{2h}^h H_{2h}^{-1}I_h^{2h} H_h \right)S^{\nu_1}
\end{equation}
where $S$ is the error iteration matrix of the smoother, and $\nu_1, \nu_2$ are the numbers of pre- and post-relaxations, respectively. 
A main challenge in the designing of multigrid methods for the Helmholtz equations lies in building a good coarse approximation $H_{2h}$.
Such an approximation should resemble the fine-grid operator $H_h$ in the spectral sense of \eqref{eq:1-lambdaRatio}.

There are two main approaches for designing a coarse approximation:
Galerkin coarsening and re-discretization. 
The Galerkin coarse grid operator is
\begin{equation}\label{eq:Galerkin}
H_{2h} = I_h^{2h}H_h I_{2h}^h.
\end{equation}
In fact, this is a projection onto the range of the interpolation:
if all the smooth error components are in the range of the interpolation, then the coarse grid correction using Galerkin coarsening annihilates the error.
High-order interpolations can represent more accurately oscillatory functions, enabling good alignment of $H_h$ and $H_{2h}$, see, e.g., \cite{dwarka2020scalable}.
However, high-order interpolation do not suffice for making $H_{4h}$ align with $H_h$ in a 3-level method.
Another drawback of the use of Galerkin coarsening with high-order interpolation is the larger effective computational stencils on the coarse grids.
We aim at modifying the Galerkin coarse grid correction so that it would be able to well-approximate the third level operator, and outperform re-discretization method despite the larger stencils.

\section{Real-shifted coarse grid correction}\label{sec:main}

Here we describe our method.
First, we present our choices of discretization and multigrid components, and then outline the application of the real shift on the third level Galerkin coarse operator.

\subsection{Discretization and multigrid components}\label{subsec:multigrid}

To discretize \eqref{eq:acousitcHelm}, we use compact 9-point (in 2D) or 19-point (in 3D) fourth-order discretizations. 
Although wider stencils can further reduce the numerical dispersion, compact stencils ease memory access and are friendly for GPU computations.
In 2D, we use the stencil 
\begin{equation} \label{eq:disc4thCompact}
\frac{1}{6h^2}
\begin{bmatrix}
-1 & -4  & -1 \\[0.4em]
-4  & 20 & -4  \\[0.4em]
-1 & -4  & -1
\end{bmatrix}
- \frac{k^2}{12}
\begin{bmatrix}
 &  1 & \\[0.4em]
 1 & 8 & 1 \\[0.4em]
 & 1 & 
\end{bmatrix},
\end{equation}
see \cite{singer1998high}, and in 3D, its analogue from \cite{turkel2013compact},
\begin{equation}
\label{eq:disc4thCompact3D}
-\frac{1}{6h^2}  
\begin{bmatrix}
\begin{bmatrix}
 & 1 & \\
 1 & 2 & 1 \\
 &  1 &
\end{bmatrix}
\begin{bmatrix}
 1 & 2 & 1  \\
2 & -24 & 2  \\
1 & 2 & 1
\end{bmatrix}
\begin{bmatrix}
 & 1 & \\
 1 & 2 & 1 \\
 &  1 &
\end{bmatrix}
\end{bmatrix}
 - \frac{k^2}{12}
\begin{bmatrix}
\begin{bmatrix}
 &  & \\
 & 1 &  \\
 &   &
\end{bmatrix}
\begin{bmatrix}
 & 1 &   \\
1 & 6 & 1  \\
  & 1 &  
\end{bmatrix}
\begin{bmatrix}
 &  & \\
 & 1 &  \\
 &   &
\end{bmatrix}
\end{bmatrix}.
\end{equation}
As coarse grid operators we use the Galerkin projections:
Let $H_h = -\Delta_h-k^2M_h$ be a discretization of \eqref{eq:acousitcHelm}.
 For the second grid, the coarse operator is given by the Galerkin coarsening \eqref{eq:Galerkin}.
The coarse grid operator  of the third grid is equipped with a real shift, as described below in Subsection \ref{subsec:realshift}.

We use 3-level methods only, for several reasons.
First, 3-level methods form a balanced compromise: for 2-level methods, the coarse grid might be too heavy, and for 4-level methods the coarse grid might not satisfy the Nyquist limit, requiring $G>2$ for any wave representation. 
Therefore, efficient dispersion correction is impossible for more than 3-levels, when starting at 10--12 grid points per wavelength on the fine grid.
Second, We aim to achieve scalability, which motivates the use of a three-level method. When the third (coarsest) level remains too large to be solved directly, it must be treated using alternative approaches that provide sufficient accuracy in order to preserve overall scalability.
It is also shown in \cite{yovel2025vanka} that if the third level is feasible to solve directly, 3-level cycles are the fastest for high-frequency Helmholtz problems.

High-order intergrid operators are shown to be efficient for Helmholtz problems \cite{dwarka2020scalable}.
We use the following stencils:
\begin{equation}\label{eq:intergrid}
I_h^{2h} = \bigotimes^d \frac{1}{16}\begin{bmatrix}
1 & 4 & 6 & 4 & 1
\end{bmatrix},
\qquad
I_{2h}^h = \frac{1}{2^d}\left(I_h^{2h}\right)^T
\end{equation}
where $d$ is the dimension and $\displaystyle \bigotimes^d$ denotes $d$ repeated applications of a Kronecker product of a stencil with itself.
In 3D, for memory reasons, we use the level-dependent intergrid scheme presented in \cite{yovel2025vanka}, which enables smaller coarse stencils.
The level-dependent intergrid method is comprised of high-order intergrid operators \eqref{eq:intergrid} between the first and second level and in the prolongation from third to second grid, while the restriction from second to third grid is of lower order:
\begin{equation}\label{eq:levdep}
I_{2h}^{4h} = \bigotimes^d \frac{1}{4}\begin{bmatrix}
1 & 2 & 1
\end{bmatrix}.
\end{equation}
Using high-order grid-transfer operators \eqref{eq:intergrid}, $H_{4h}$ and $M_{4h}$ are represented by $7\times 7$ points stencils (in 2D) or $7\times 7\times 7$ points stencils (in 3D).
Similarly, using level-dependent intergrid operators yields $5\times 5$ and $5\times 5\times 5$ points stencils, respectively.

\subsection{Application of the real shift}\label{subsec:realshift}

Here we describe how to apply the real shift to the third grid.
The high-order discretizations \eqref{eq:disc4thCompact} and \eqref{eq:disc4thCompact3D}, combined with Galerkin coarsening \eqref{eq:Galerkin} using the high-order intergrid operators \eqref{eq:intergrid}, already yield well-aligned phases of the first and the second grid.
To align the third grid as well, we impose a real shift in the following manner. 

With $H_h = -\Delta_h-k^2M_h$ as before, let $\alpha$ be a correcting factor to be optimized later.
We define the real-shifted fine grid operator
\begin{equation}\label{eq:realShift}
H_h^\alpha = -\Delta_h-(\alpha k)^2M_h = H_h +(1-\alpha^2)k^2M_h
\end{equation}
and then apply the Galerkin coarsening twice to obtain the real-shifted coarse grid operator:
\begin{equation}\label{eq:CoarseRealShifted}
H_{4h}^\alpha = I_{2h}^{4h}\left(I_h^{2h}H_h^\alpha I_{2h}^h\right)I_{4h}^{2h}.
\end{equation}
We note that the shift parameter $\alpha$ in $H_{4h}^\alpha$ does something more sophisticated than simply stretching or translating the eigenvalues on the real axis.
In fact, for constant $k$,
\begin{equation}\label{eq:realShiftMeaning}
H_{4h}^\alpha = H_{4h} + (1-\alpha^2)k^2 M_{4h}
\end{equation}
where $H_{4h}$ and $M_{4h}$ are the Helmholtz and mass operators after undergoing two Galerkin coarsenings, resulting in wider stencils.
Our method exploits the advantages of wide stencils, while keeping the fine grid discretization compact.

In the next section we present analysis that enables to choose $\alpha$ such that the relative dispersion error, or the relative maximal phase difference between $H_h$  and $H_{4h}^\alpha$, is minimized over  all possible propagation directions.
The choice of $\alpha$ depends only on $G$, the number of grid points per wavelength on the fine grid, and hence $\alpha$ can be optimized once and then the same parameter is used for all the grid sizes (while keeping $G$ constant).

\section{Grid-to-grid dispersion analysis}\label{sec:analysis}

We present a dispersion analysis method tailored to our setting. 
In classical dispersion analysis, the dispersion error, which is defined to be the phase difference between the continuous operator and its discretization, is analyzed in Fourier space.
Unlike classical dispersion analysis, we compare grids directly rather than comparing to the continuous case.
That is, we define a new sense of dispersion error, comparing the first and the third grid discrete operators.
This requires grid stretches since the different grids has different corresponding numbers of points per wavelength.
Unlike classical three-grid local Fourier analysis, our analysis does not include the smoother and intergrid application.
This isolates the dispersion effect and enables finding the correct real shift for a stable coarse grid correction, rather than predicting overall convergence.
That is, our method cannot be used for predicting convergence rate, but can be used for predicting the optimal shift and predict when, no matter how good the smoother will be, the method will fail.
We give a rule of thumb to analyze when the pollution effect poses a limitation.
For every given error, there exists a grid size for which the error accumulates and the coarse grid correction  fails.
We asses this limiting grid size.

In Subsection \ref{subsec:disp_back} we give essential definitions from classical dispersion analysis. 
In subsection \ref{subsec:grid2grid} we present our grid-to-grid analysis, and in Subsection \ref{subsec:examples} we show numerical examples of the resulting dispersion relation. 
In particular, we use the analysis to tune the real-shift to be used for various scenarios and predict their limitations. 

\subsection{Classical dispersion analysis}\label{subsec:disp_back}
We start with some background on dispersion analysis.
The \emph{dispersion relation} of a continuous or discrete operator $H$ is the set of frequencies for which the symbol vanishes. 
For dimension $d$, it takes the form
\begin{equation}\label{eq:disp_rel}
\left\{\theta\in[-\pi,\pi]^d \mid \widetilde{H}(\theta) = 0\right\}
\end{equation}
where $\widetilde{H}(\theta)$ is the symbol of $H$.

The dispersion relation of the continuous Helmholtz equation is always a circle (in 2D) or a sphere (in 3D) because of the spherical symmetry of the Laplacian.
In more detail, the dispersion relation of the continuous Helmholtz equation with wavenumber $k$ is a circle of radius $k$:
\begin{equation}\label{eq:disp_rel_cont}
\{\theta\in\mathbb R^d \mid |\theta|^2 - k^2 = 0 \} = \{\theta\in\mathbb R^d \mid |\theta| = k \}
\end{equation}
see, e.g., \cite{cocquet2021closed}.
When we compare it to the discrete dispersion relation, $h$ need to be taken into account.
We work with the dimensionless wavenumber $kh$ and then the radius of the circle is equal to $r = kh = 2\pi/G$, where $G$ is the number of grid points per wavelength.

For a discrete matrix given as a convolution against a stencil, the symbol can be calculated as a sum of complex exponents. 
For 2D compact stencils:
\begin{equation}\label{eq:general_stencil}
S = \begin{bmatrix}
s_{-1,1} & s_{0,1} & s_{1,1} \\
s_{-1,0} & s_{0,0} & s_{1,0} \\
s_{-1,-1} & s_{0,-1} & s_{1,-1} 
\end{bmatrix}
\quad \Rightarrow \quad
\widetilde{H}_h(\theta) = \sum s_{k,l} e^{\im (k\theta_1+l\theta_2)},
\end{equation}
see, e.g., \cite{cocquet2021closed}.
Applying \eqref{eq:general_stencil} to the stencil \eqref{eq:disc4thCompact} gives
\begin{equation}\label{eq:symbol_fourth}
\widetilde{H}_h(\theta) =
\frac{1}{h^2}\left(
 \frac{10}{3}-\frac{2}{3}(kh)^2 - \left(\frac{4}{3}+\frac{1}{6}(kh)^2\right)\left(\cos(\theta_1)+\cos(\theta_2)\right) -\frac{2}{3} \cos(\theta_1)\cos(\theta_2)
 \right)
\end{equation}
and similarly in 3D, the same procedure yields
\begin{equation}\label{eq:symbol_fourth_3D}
\widetilde{H}_h(\theta) =
\frac{1}{h^2}\left(
4-\frac{1}{2}(kh)^2 - \left(\frac{2}{3}+\frac{1}{6}(kh)^2\right)\sum_{i=1}^3\cos(\theta_i) -\frac{2}{3}
 \sum_{1\leq i < j\leq 3} \cos(\theta_i)\cos(\theta_j)
 \right)
\end{equation}
when applied to the stencil \eqref{eq:disc4thCompact3D}.
The discrete dispersion relation is not a perfect circle or sphere, since $H_h(\theta)$ from \eqref{eq:symbol_fourth} or \eqref{eq:symbol_fourth_3D}, 
is not a radial function.
Hence, numerical dispersion can never be fully corrected by a shift in more than 1D \cite{ernst2013multigrid}.

To compute the discrete radius $r_1(\phi)$ of the fine grid operator when $\phi$ is a propagation direction, we sample the symbol over a ray of direction $\phi$ in $[-\pi,\pi]^d$.
For instance, in 2D we sample $\theta\in[-\pi,\pi]^2$ for which $\arctan(\theta_2/\theta_1)=\phi$. 
In 3D, we refer to $\phi$ as a vector, whose entries are azimuhtal and polar angles, and we sample $\theta$ accordingly.
Subsequently, we define the discrete radius $r_1(\phi)$ as the distance between the origin and the point at which the symbol switches sign.
The ratio of radii is the inverse ratio of phase velocities
\begin{equation}\label{eq:inverse_ratio}
\frac{r}{r_1(\phi)} = \frac{V_1(\phi)}{V}
\end{equation}
where $V$ be the continuous wave speed and $V_1(\phi)$ the direction-dependent numerical wave speed of the fine-grid operator.

The \emph{dispersion error} is defined (see, e.g., \cite{jo1996optimal}) as the normalized phase difference
\begin{equation}\label{eq:disp_err}
e_{d}(\phi) = \frac{V_1(\phi) - V}{V} = \frac{V_1(\phi)}{V}-1
\end{equation}
where the velocity ratio can be calculated by the radii ratio through \eqref{eq:inverse_ratio}.

\subsection{Our grid-to-grid analysis}\label{subsec:grid2grid}
In our context, we want to measure how well-aligned are the first and third grids in a 3-level method.
It suffices for predicting optimal convergence since the first two grids are already aligned when using high-order intergrid operators \cite{dwarka2020scalable}.
Hence, we define the following measure of error:
\begin{equation}\label{eq:disp_grids_err}
e_{g}(\alpha,\phi) = \frac{V_1(\phi)}{4V_3(\alpha,\phi)}-1
\end{equation}
where $V_3(\alpha,\phi)$ is the numerical wave speed of the third-grid Galerkin operator $H_{4h}^\alpha$ from \eqref{eq:CoarseRealShifted} at propagation direction $\phi$.

The factor 4 in \eqref{eq:disp_grids_err} compensates the two coarsening steps by which $V_1$ and $V_3$ differ:
since the coarse grid operator has $G/4$ grid points per wavelength, its dispersion relation is a curve with direction-dependent radius which is approximately 4 times larger compared to that of the fine grid.
We want to see the alignment after stretching the fine radius by 4, since this represents the quality of the coarse approximation.
Multiplying the direction-dependent radius of the fine grid by 4 is equivalent to dividing the corresponding phase velocity $V_1(\phi)$ by 4, since the velocities ratio is the inverse of the radii ratio.

Next, we explain how to calculate $e_g(\alpha,\phi)$ from \eqref{eq:disp_grids_err}.
To calculate $V_1(\phi)$  we use \eqref{eq:general_stencil}, and to calculate $V_3(\alpha,\phi)$ we first extract the effective stencil that corresponds to the Galerkin operator \eqref{eq:CoarseRealShifted}. 
This is done by stencil-to-stencil convolutions and down-sampling, as shown in \cite{chen2024matrix}.
That is, we first calculate the wide stencil corresponding to $H^\alpha_{2h}$ by convolving the stencils of $I_h^{2h}$, $H^\alpha_h$ and $I_{2h}^h$ (from an $\alpha$-shifted version of  \eqref{eq:Galerkin}) and then we derive the stencil corresponding to $H^\alpha_{4h}$ by convonving the relevant  intergrid stencils against the resulting $H^\alpha_{2h}$.
Once we have the stencil of the Galerkin coarse operator $H^\alpha_{4h}$, we calculate the symbol by a wider version of \eqref{eq:general_stencil}.

The accumulation of the error $e_g$ over wavelengths forms a main limitation to any dispersion correction method.
We give here heuristic quantitative bounds.
Let $G$ be the number of grid points per wavelength and let $n_{crit}$ be the maximal number of grid points per direction in 2D or in 3D for which the method converges. 
Let $e_g=e_g(\alpha)$ be the error from \eqref{eq:disp_grids_err}. 
Starting the accumulation from a peak (the wave source), the accumulated error over $l$ wavelengths is $le_g$.
Reasonable assumptions\footnote{After less then 1/4 wavelength, the wave haven't change sign, assuming a peak in the source's location. 
After 1/2 wavelength, the wave reached an opposite phase which exposes it to a maximal error. 
Thus, it is reasonable to assume that the relative part of a wavelength for which there is enough error to interfere with convergence is between these two numbers.} lead us to demand that the accumulated error is roughly between $1/4$ and $1/2$ wavelength at $n_{crit}$, and then substituting $G=n_{crit}/l$ yields the bound
\begin{equation}\label{eq:n_crit_bound}
\frac{G}{4e_g} \leq n_{crit} \leq \frac{G}{2e_g}.
\end{equation}
This bound can be thought of as a rule of thumb, to assess the largest grids that can be solved for a certain theoretical error. 
We show in the next subsection that these bounds are obtained numerically.
Improving them is a subject for our future research.

\subsection{Numerical examples}\label{subsec:examples}

\begin{figure}
\begin{center}
	\newcommand{\image}[1]{\includegraphics[width=0.49\linewidth]{./#1}}
    \subfigure[\footnotesize $G=12$, $\alpha=0$]{\image{disp_12G_nonshifted.eps}\label{fig:disp_12G_nonshifted}}
    \subfigure[\footnotesize $G=12$, $\alpha=\alpha^*$]{\image{disp_12G_shifted.eps}\label{fig:disp_12G_shifted}}
\end{center}
\caption{Dashed curve: dispersion relation \eqref{eq:disp_rel} of the coarse grid operator \eqref{eq:CoarseRealShifted} for the 2D Helmholtz equation with 12 grid points per wavelength on the fine grid and level-dependent intergrid.
Solid curve: dispersion relation of the corresponding fine grid operator stretched by a factor of 4.
}\label{fig:disp_rel}
\end{figure}

\begin{figure}
  \begin{center}
    \includegraphics[width=0.7\textwidth, trim=2cm 8cm 2cm 8cm, clip]{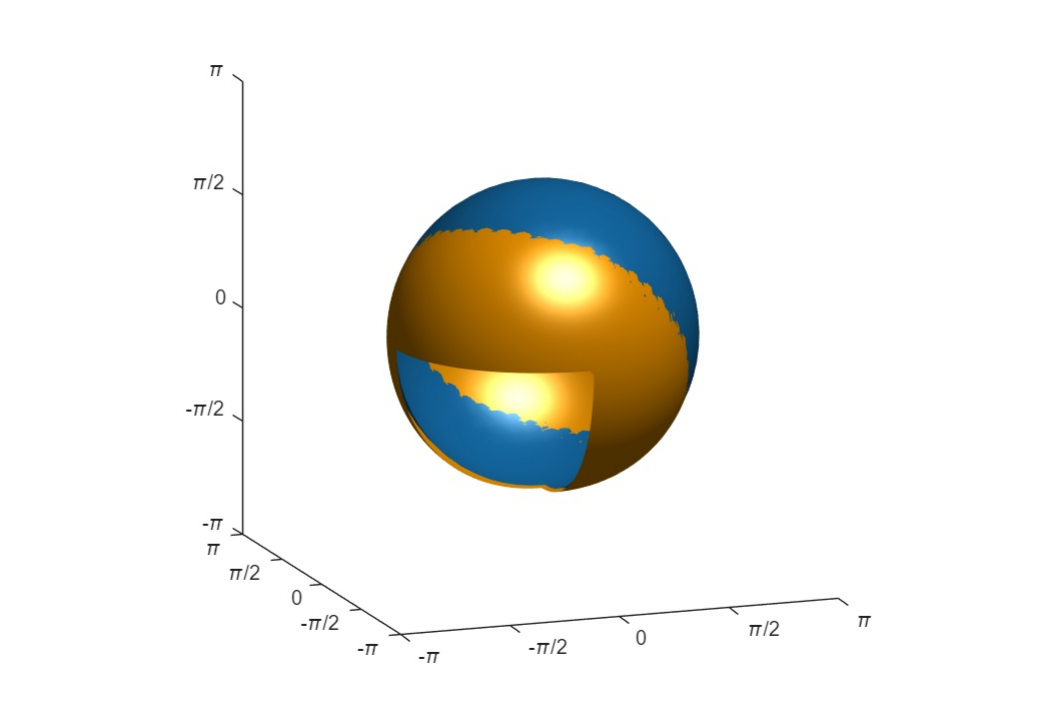}\\
  \caption{Orange: dispersion relation \eqref{eq:disp_rel} of the coarse grid operator \eqref{eq:CoarseRealShifted} for the 3D Helmholtz equation with 12 grid points per wavelength on the fine grid and level-dependent intergrid.
Blue: dispersion relation of the corresponding fine grid operator stretched by a factor of 4.}\label{fig:disp_rel_3D}
  \end{center}
\end{figure}

We define the optimal real shift $\alpha$ as the solution of the following min-max problem:
\begin{equation}\label{eq:minmax}
\arg\min_{\alpha}\max_{\phi} (e_{g}(\alpha,\phi))
\end{equation}
which is equivalent to minimizing the infinity norm of the vector of errors. 
We solve \eqref{eq:minmax} by exhaustive search on $\alpha$. 
One can consider using more efficient numerical optimization methods, however, brute-force is feasible since we optimize only one parameter. 
The sought parameter $\alpha^*$ can be then used to solve problems of different grid sizes, domains and media, as demonstrated in Section \ref{sec:numerical}.

To find an optimal $\alpha$ in 2D, we sample \eqref{eq:minmax} for $\phi\in[0,\pi/4]$, with a resolution of $0.1$. 
Similarly, in 3D, we sample the azimuthal angle $\phi_1 \in [0,\pi/4]$ and the polar angle $\phi_2\in[\pi/2-\arccos(1/\sqrt{3}),\pi/2]$, each with a resolution of $0.1$.
The remaining propagation directions produce the same errors due to symmetry.
For each propagation direction we sample $\alpha$ with a resolution of $5\cdot10^{-4}$, and for each $\phi$ and $\alpha$ we sample the symbol on the ray of direction $\phi$ in the $\theta$ plane with a resolution of $10^{-3}$. 
The sampling in the $\theta$ plane reads, in 2D, sampling for $\theta\in[-\pi,\pi]^2$ satisfying $\theta_1=r\cos(\phi)$ and $\theta_2=r\sin(\phi)$ for some $r$ and for the same $\phi$.

\begin{table}
\centering
\begin{tabular}{c|ccc|ccc}
  \toprule
  \mc{7}{c}{Optimal real shift and error in 2D} \\
  \midrule
 Intergrid & \mc{3}{c|}{ cubic} & \mc{3}{c}{ lev-dep} \\
  \midrule
&  $G=10$ & $G=11$ & $G=12$ &  $G=10$ &  $G=11$ &  $G=12$ \\
  \midrule
  $\alpha^*$ &  1.0140 & 1.0075 & 1.0045 & 1.0290 & 1.0190 &  1.0135 \\
 $ \max_{\phi}e_{g}^*\times10^{2}$ &  1.1924 & 0.6130 &  0.3340 &  1.7117 &  1.1821 & 0.8111 \\
  $n_{crit}\in$  &  [209,419]  &  [449,897] &  [898,1796] & [146,292] & [151,302] & [369,739]  \\
$n_{crit}$ num &  248 & 536 &  1396 & 184 &  284 & 508  \\
 \bottomrule
 \end{tabular}
\caption{Optimal solution and objective for \eqref{eq:minmax}, and the resulting bounds on the  critical grid size for different numbers of grid points per wavelength in 2D.
The numerical $n_{crit}$ is the largest grid size for which the multigrid convergence is stable for a point source located at the center of the upper boundary.
}
\label{tab:disp_err}
\end{table}

\begin{table}
\centering
\begin{tabular}{c|ccc|ccc}
  \toprule
  \mc{7}{c}{Optimal real shift and error in 3D} \\
  \midrule
Intergrid & \mc{3}{c|}{ cubic} & \mc{3}{c}{ lev-dep} \\
  \midrule
&  $G=10$ & $G=11$ & $G=12$ & $G=10$ & $G=11$ & $G=12$ \\
  \midrule
$\alpha^*$ & 1.0130 & 1.0065 & 1.0045 & 1.0245 & 1.0165 &  1.0120  \\
$ \max_{\phi}e_{g}^*\times10^{2}$ & 1.2739 & 0.6649 & 0.3340 &  2.0668 & 1.3369 &  0.9542 \\
$n_{crit}\in$ & [196,392] & [414,827] & [898,1796] & [120,241] & [201,402] &  [314,628]  \\
 \bottomrule
 \end{tabular}
\caption{Optimal solution and objective for \eqref{eq:minmax}, and the resulting bounds on the  critical grid size for different numbers of grid points per wavelength in 3D.
}
\label{tab:disp_err_3D}
\end{table}

In Table \ref{tab:disp_err} we report the optimal real shifts and the resulting errors for different values of $G$ in 2D.
We observe that the maximal error for optimal $\alpha$ is of order of $10^{-2}$ for $G=10$ and of order of $10^{-3}$ for $G=12$, when using cubic intergrid.
This difference, in fact, makes it impossible to retain scalability for large grids using our method with 10 grid points per wavelength: 
the anisotropy of the resulting Galerkin stencil is large, and hence the dispersion can be only corrected to some extent by a real shift, and the accumulation of phase shift interferes with convergence even for relatively small grids.
To assess this phenomenon quantitatively, we report the theoretical bounds on the critical number of grid points per direction $n_{crit}$.
We show that the numerical values of $n_{crit}$ are within these bounds.
Note that $n_{crit}$ and its bounds were calculated assuming a wave that travels the entire length of the domain. 
If the point source is located in the middle of the domain, we observe in practice a larger $n_{crit}$.

In Table \ref{tab:disp_err_3D} we report the optimal shifts, errors and expected bounds on $n_{crit}$ in 3D.
Although the results are very similar to the 2D table \ref{tab:disp_err}, it is worth mentioning that the same limitations on $n_{crit}$ value can be less significant in 3D, where memory considerations might restrict the user of the method to smaller number of grid points per direction.

In Fig. \ref{fig:disp_rel}
we depict the 2D dispersion relations with and without the shifting, corresponding to the rightmost column of Table \ref{tab:disp_err}. 
The sharp-eyed reader can see the slight misalignment in Fig. \ref{fig:disp_12G_nonshifted} which is corrected in \ref{fig:disp_12G_shifted}. 
Similarly, Fig. \ref{fig:disp_rel_3D} shows the 3D dispersion relations corresponding to the rightmost column of Table \ref{tab:disp_err_3D}, showing good alignment of the coarse grid with the stretched dispersion relation of the fine grid.

\begin{figure}
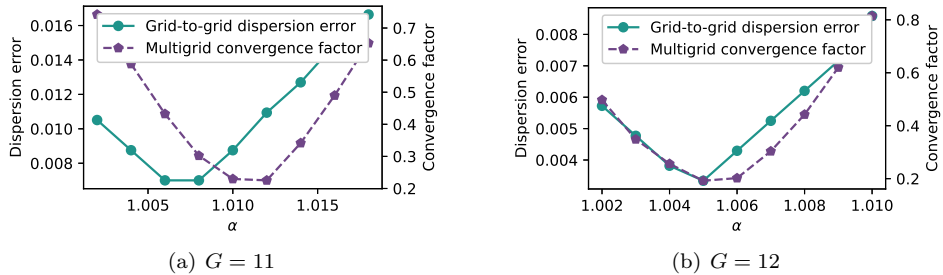

\begin{center}
	\newcommand{\image}[1]{\includegraphics[width=0.47\linewidth]{./#1}}
    \subfigure[\footnotesize $G=11$]{\image{opt_shift_11G.eps}\label{fig:opt_shift_11G}}
    \hspace{10pt}
    \subfigure[\footnotesize $G=12$]{\image{opt_shift_12G.eps}\label{fig:opt_shift_12G}}
\end{center}
\caption{Multigrid convergence factor and grid-to-grid dispersion error vs. the real shift $\alpha$ from \eqref{eq:CoarseRealShifted}, for different $G$ values. 
The convergence factor was calculated for a grid size of $256^2$ cells for $G=10$, $512^2$ cells for $G=11$ and $1024^2$ cells for $G=12$.
}\label{fig:opt_shift}
\end{figure}

Finally, we investigate the alignment of the theoretically optimal shift and the numerically optimal shift.
We compute the grid-to-grid dispersion error as a function of the shift $\alpha$ and compare it to the multigrid convergence factor, measured for a grid size which is slightly lower than $n_{crit}$.
The results are displayed in Fig. \ref{fig:opt_shift}. 
In Fig. \ref{fig:opt_shift_12G} we see great alignment between the prediction of the optimal shift and the convergence in practice for $G=12$ and a grid size of $1024^2$ cells.
For $G=11$ and a grid of $512^2$ cells, the alignment slightly deteriorates, see Fig. \ref{fig:opt_shift_11G}.
The theory predicts an optimal shift of $\alpha = 1.008$ whereas the optimal shift in practice is $\alpha=1.012$.
This less accurate prediction might occur because of the smaller grid size. 
The theoretical predictions in the Fourier space assume an infinite domain (or periodic boundary conditions) and hence are generally more accurate for large grids. 
However, because of the smaller $n_{crit}$, we could not achieve a much larger grid size and while retaining convergence over the examined range of shifts.

\begin{remark}
One could consider, for memory reasons, a similar method with bilinear intergrid. 
Such a method would require introducing real shifts for both coarse grids, and the resulting errors after the shifts are still too large.
For $G=12$, e.g., the optimal error is 
$\max_{\phi}e_{g}^* = 0.0361$ and the resulting expected critical grid size is $83\leq n_{crit} \leq 166$ which is too small for most reasonable uses.
Hence, we consider the method as not applicable for standard intergrid operators.
\end{remark}

\section{Numerical results}\label{sec:numerical}

In this section we give numerical experiments to demonstrate the efficiency of our method. 
Throughout the experiments, we solve the Helmholtz equation equipped with a padding layer of 20 cells of ABC. 
We use convergence tolerance of relative residual $<10^{-6}$.
Our code is written in \texttt{Julia} language \cite{Julia}.
The experiments are calculated on a dual-core laptop with 32 GB RAM running Windows 11.

\subsection{Two-dimensional experiments}\label{subsec:2d}

The multigrid framework for the 2D experiments includes 3-level $W(1,1)$ cycles with damped Jacobi smoother, which is cheap, parallel and easy to implement.
We use the optimal damping parameter reported in \cite{yovel2025vanka}, which is $0.89$ for both the first and the second grid relaxation.
Only one pre- and one post-relaxation are applied, since the smoother is non-convergent for a larger number of relaxations in the absence of a complex shift.
Unless stated otherwise, we employ high-order intergrid operators, see \eqref{eq:intergrid}, that are shown to reduce numerical dispersion (see Section \ref{sec:main}) and run W-cycles which better preserve scalability.
As we show later, however, V-cycles are nearly as efficient in our setting and may even outperform them for small grids.

As coarse grid operators, the modified coarse approximation developed in Section \ref{sec:main} is used.
That is, the Galerkin projection \eqref{eq:Galerkin} is used for the second level and \eqref{eq:CoarseRealShifted} for the third (and coarsest) level, with the shifts $\alpha$ found by our theory as reported in Table \ref{tab:disp_err}.
We refer to our method as RS-CGC (real-shifted coarse grid correction) and compare it the classical CSLP (complex-shifted Laplacian preconditioner) from \cite{erlangga2006novel}.

\subsubsection{Homogeneous 2D media}

\begin{figure}
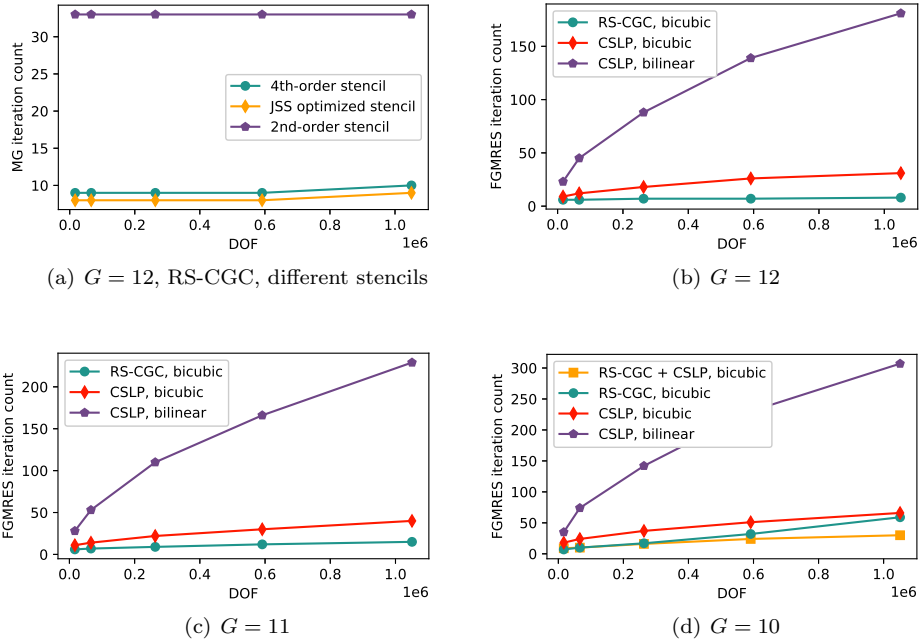

\begin{center}
	\newcommand{\image}[1]{\includegraphics[width=0.49\linewidth]{./#1}}	\subfigure[\footnotesize $G=12$, RS-CGC, different stencils]{\image{disc_comp.eps}\label{fig:disc_comp}}
    \subfigure[\footnotesize $G=12$]{\image{cslp_rsms_comp.eps}\label{fig:cslp_rsms_comp12}}
\subfigure[\footnotesize $G=11$]{\image{cslp_rsms_comp_11G.eps}\label{fig:cslp_rsms_comp11}} 
    \subfigure[\footnotesize $G=10$]{\image{cslp_rsms_comp_10G.eps}\label{fig:cslp_rsms_comp10}}
\end{center}
\caption{FGMRES(20) iteration count for the homogeneous 2D Helmholtz equation, where a 3-level W(1,1) cycle with Jacobi smoothing serves as a preconditioner, for CSLP, RS-CGC and their combination.
The complex shifts for the CSLP are: for bilinear intergrid operators, $0.3k^2$ and for bicubic intergrid operators, $0.1k^2$. For RS-CGC $+$ CSLP, the complex shift is $0.03k^2$.
}\label{fig:cslp_rsms_comp}
\end{figure}

In the first experiment, we compare the performance of our method with different stencils used for the discretization of the 2D Helmholtz equation \eqref{eq:acousitcHelm}.
We use the standard second-order stencil, the fourth-order stencil \eqref{eq:disc4thCompact} and the dispersion-correction stencil presented in \cite{jo1996optimal},
\begin{equation}\label{eq:JSS_sten}
H_{JSS} = -\frac{1}{h^2}\left(a \begin{bmatrix}
 & 1 & \\
1 & -4 & 1 \\
 & 1 &
\end{bmatrix}
+ \frac{1-a}{2}
\begin{bmatrix}
1 &  & 1\\
 & -4 &  \\
1 &  & 1
\end{bmatrix}\right)
- k^2 \begin{bmatrix}
\frac{1-b-c}{4} & \frac{c}{4} & \frac{1-b-c}{4} \\[0.5em]
\frac{c}{4} & b & \frac{c}{4} \\[0.5em]
\frac{1-b-c}{4} & \frac{c}{4} & \frac{1-b-c}{4} 
\end{bmatrix}
\end{equation}
with $[a,b,c] = [0.5461,0.6248,0.37524]$. 
These weights were reported in \cite{jo1996optimal} to be optimal in terms of $L^2$ norm of the dispersion error, summed over all possible propagation directions, for $G>4$.
We refer to this as the JSS stencil.
We note that the JSS stencil is only second-order accurate.
The results are depicted in Fig. \ref{fig:disc_comp}.
Interestingly, all the stencils keep the same scalability behavior, even the standard second-order stencil which is characterized by high numerical dispersion. 
We explain this by the behavior of the Galerkin coarsening, characterized by wider stencil on the corase grids, and good alignment with the fine-grid operator.
When the the real shift is interwoven into the Galerkin operator, it forms a stable dispersion correction even for basic stencils.
However, the iteration count for the standard second-order stencil is significantly larger.
The fourth-order stencil and the JSS stencil gives similar performance.
Despite the slight improvement that JSS offers, in the remainder of this paper we use the fourth order discretization \eqref{eq:disc4thCompact} in 2D and \eqref{eq:disc4thCompact3D} in 3D, which gives a lower order of truncation error in general, and 3D generalization that does not require re-optimizing weights.

\begin{remark}
In \cite{jo1996optimal}, the objective of the optimization process was the squared error integrated over propagation directions and $G$ values. 
There, the discretization was developed for general purposes, and the minimization was performed in an $L^2$ sense, which represents one possible choice among many.
In the context of multigrid, however, it may be more appropriate to control the error in an $L^\infty$ sense, since even a single propagation direction with a large error may deteriorate the multigrid convergence.
\end{remark}

In the second experiment, we compare our method to the well-known CSLP \cite{erlangga2006novel} with standard and with high-order grid-transfer operators.
For the CSLP frameworks, we use small complex shifts that optimize convergence:
the shifts $0.1k^2$ for high-order intergrid operators, and $0.3k^2$ for standard intergrid operators, were chosen by trial and error.
Fig. \ref{fig:cslp_rsms_comp12} demonstrates perfect scalability of RS-CGC for 12 grid points per wavelength on the fine grid.
Our method converges in 6--8 iterations --- wavenumber independent convergence --- up to a grid resolution of $1024^2$ cells.
In Fig. \ref{fig:cslp_rsms_comp11} we repeat the experiment for $G=11$. 
Here, RS-CGC converges in 6--15 iterations for the different grids, $2.6$ to $15.2$ times faster than CSLP with high-order and standard intergrid operators, respectively.
However, for $G=10$ the scalability breaks at $512^2$ cells, see Fig. \ref{fig:cslp_rsms_comp10}. 
The RS-CGC method remains $5$ times faster than classical CSLP and equally efficient as CSLP with high-order intergrid.
We suggest combining RS-CGC with CSLP, which achieves a speedup by an order of magnitude compared to classical CSLP.
The optimal complex shift for the combination of RS-CGC with CSLP, determined by trial and error, is only $0.03k^2$. This may explain the improved scalability: the larger the shift, the farther the preconditioner is from the original operator.

\begin{figure}
  \begin{center}
    \includegraphics[width=0.45\textwidth]{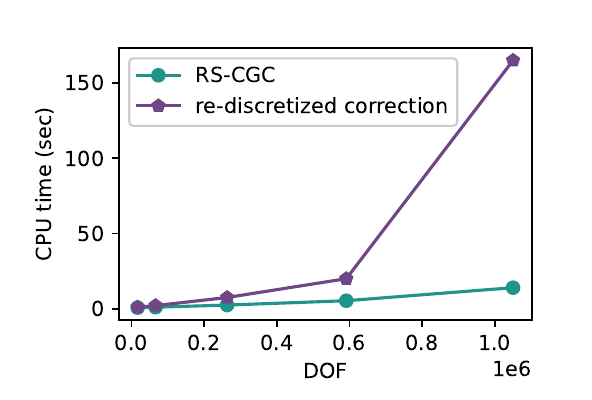}\\
  \caption{
  Time vs. DOF for solving the 2D Helmholtz equation \eqref{eq:acousitcHelm} with $G=16$ using preconditioned GMRES(20) with RS-CGC or a re-discreitization cycle as a preconditioner.
  }
  \label{fig:re_disc_comp}
  \end{center}
\end{figure}

In the third experiment, we compare RS-CGC to a re-discretization based dispersion correction method.
The re-discretization method is inspired by \cite{cocquet2017finite}, where the reported optimal values for a stencil of the form \eqref{eq:JSS_sten} with $G=4$ are $[a,b,c] = [0.6054,1.0532,0.0002],$ and the optimal modified wavenumber is $k^* = 0.87725k$.
We solve the 2D Helmholtz equation \eqref{eq:acousitcHelm} for $G=16$ on the fine grid (and hence $G=4$ on the coarse grid), for different grid sizes, using preconditioned GMRES(20). 
For the re-discretized method we use bilinear intergrid operators and re-discretization with the fourth-order stencil \eqref{eq:disc4thCompact} on the first and second level. 
On the third level we use the stencil \eqref{eq:JSS_sten} and the above mentioned weights and modified wavenumber.
It can be seen as an analogue of the method described in \cite{cocquet2021closed}, only that there 10--20 cycles of Kaczmarz relaxation were applied at each level of the cycle.
We compare this method, comprised of compact stencils over all of the levels and standard intergrid operators, to our RS-CGC. 
The aim of the comparison is to see if the high-order interpolations and wider Galerkin stencils that we use ``pay off'' in terms of time measurements.
In Fig. \ref{fig:re_disc_comp} we see that despite the Galerkin coarsening resulting with wider stencils, our method is faster and scales better. 
We only examined relatively high $G$ values, since the re-discretized method did not converge for lower $G$ values.
Even in this regime, our method scales better, and  is $\times 3.75$ faster for relatively small grids and $\times 11.8$ faster for the largest examined grid.

\subsubsection{Heterogeneous 2D media}

\begin{figure}
\begin{center}
	\newcommand{\image}[1]{\includegraphics[width=0.3\linewidth, trim=5.5cm 8cm 4cm 8cm, clip]{./#1}}
    \subfigure[\footnotesize Wedge model]{\image{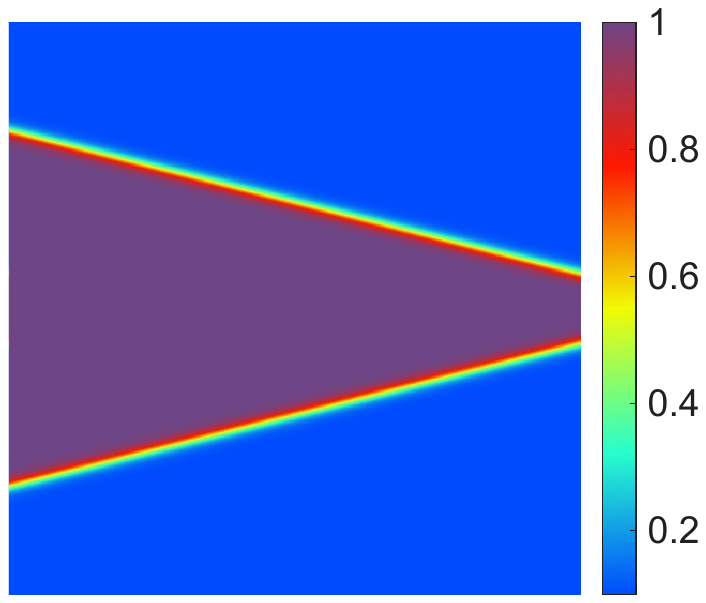}\label{fig:wedge}}
    \hspace{50pt}
    \subfigure[\footnotesize Linear model]{\image{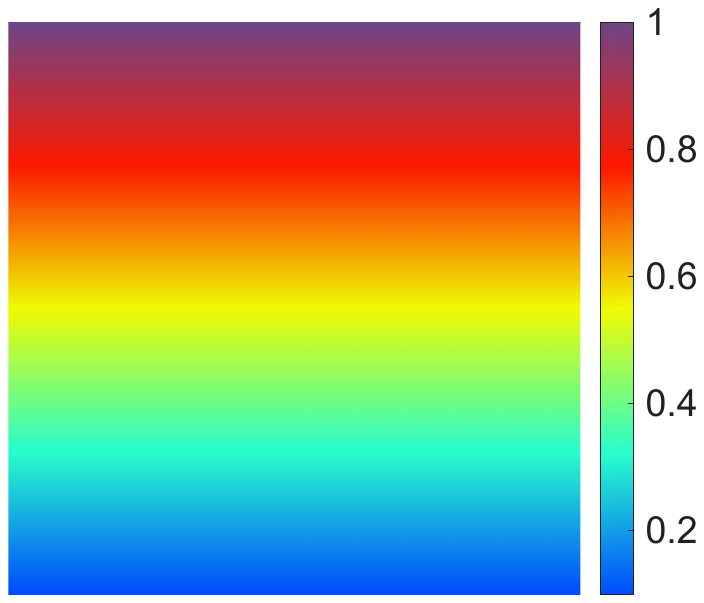}\label{fig:linear}}\\
\end{center}
\caption{On the left, the wedge velocity model and on the right, the linear velocity model on the dimensionless unit square, with slowness squared $\kappa^2\in[0.1,1]$.
}\label{fig:wedge_linear}
\end{figure}

First, we show that similar scalability holds for heterogeneous media.
We test velocity models with linearly varying $\kappa$, as well as the wedge velocity model shown in Fig.~\ref{fig:wedge_linear}. 
For both models, we consider the ranges $\kappa \in [0.5,1]$ and $\kappa \in [0.1,1]$.
Table \ref{tab:2d} demonstrates that our method is robust to different variations in coefficients.
In heterogeneous models, ensuring $G=12$ in the most difficult region implies higher $G$ values in the rest of the domain. 
Therefore, these problems are, in practice, less challenging than the homogeneous model problem.
We note that our method uses a multigrid cycle applied to $H$ (and not to a shifted version), hence it can also be used as a multigrid solver, and not only preconditioner.
The numbers in parentheses in Table \ref{tab:2d} show that the convergence of the multigrid is scalable and similar to the convergence of the GMRES-preconditioned method.

\begin{table}
\centering
\begin{tabular}{c|cccc}
  \toprule
  \mc{5}{c}{Iteration count for 2D heterogeneous media, $G=12$}\\
 \midrule
 Grid size (cells) &  linear $[0.5,1]$ &  linear $[0.1,1]$ &  wedge $[0.5,1]$  &  wedge $[0.1,1]$  \\
  \midrule
 $128\times128$  & 6 (9) & 6 (8) & 6 (9) & 6 (9) \\
 $256\times256$  & 6 (9) & 6 (8) & 6 (9) & 6 (9) \\
 $512\times512$  & 6 (9) & 6 (8) & 7 (9) & 7 (9) \\
 $1024\times1024$  & 7 (9) & 8 (9) & 8 (9) & 9 (10) \\
  \bottomrule
 \end{tabular}
\caption{Non-restarted FGMRES iteration count for solving the 2D Helmholtz for linear and wedge models of different slowness ranges.
The RS-CGC 3-grid W(1,1) cycle serves as a preconditioner.
In parentheses: iteration count for RS-CGC as a solver.
}
\label{tab:2d}
\end{table}

\begin{figure}
  \begin{center}
    \includegraphics[width=0.6\textwidth, trim=2cm 11cm 2cm 11cm,clip]{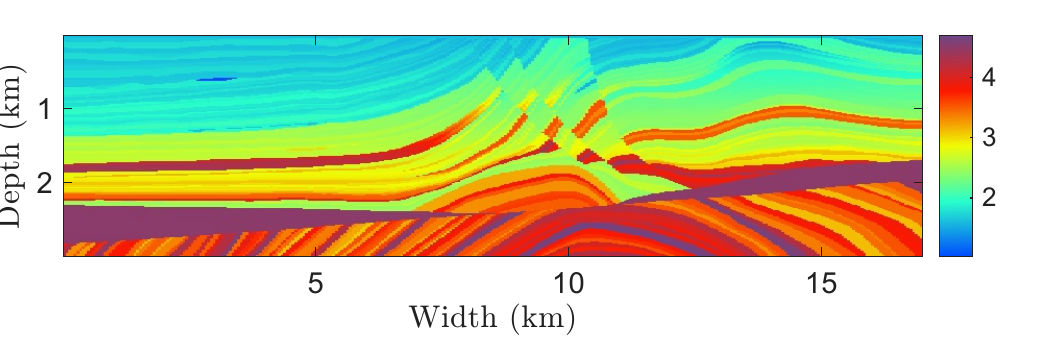}\\
  \caption{Marmousi model. Velocity units: $km/sec$.}\label{fig:Marmousi}
  \end{center}
\end{figure}

Second, we test our method on the Marmousi model \cite{brougois1990marmousi}. 
This model is 2D but resembles real world data, as it is based on a section of the Kwanza basin in Angola. 
The original domain is very shallow: its size is $17\times3.5km$.
Hence, we extend the model downwards by 12--16 cells to enable the application of the absorbing boundary layer without adding large artificial attenuation into the domain.
In Table \ref{tab:Marmousi} we compare the performance of RS-CGC for V-cycles and W-cycles. 
For small grids, the methods perform similarly in terms of iteration count, and hence the V-cycles are faster.
However, for the largest examined grid the scalability of the method deteriorates when using V-cycles.

%

\begin{table}
\centering
\begin{tabular}{c|cc|cc|cc|cc}
  \toprule
  \mc{9}{c}{Iteration count for 2D Marmousi media}\\
 \midrule
 & \mc{4}{c|}{$G=12$, RS-CGC} & \mc{4}{c}{$G=10$, RS-CGC + CSLP} \\
  \midrule
  &  \mc{2}{c|}{W-cycles} &  \mc{2}{c|}{V-cycles}  &  \mc{2}{c|}{W-cycles} &  \mc{2}{c}{V-cycles} \\
 \small Grid size (cells) & \small \#Iter & \small time (s) &  \small \#Iter & \small time (s) & \small \#Iter & \small time (s) &  \small \#Iter & \small time (s) \\
  \midrule
 $544\times128$  & 7 & 0.894 & 8 & 0.773 & 8 & 0.968 & 9 & 0.930 \\
 $1088\times240$  & 7 & 2.056  & 8 & 1.746 & 10 & 3.021 & 15 & 3.275 \\
 $2176\times464$  & 9 & 10.365 & 11 & 8.757 & 19 & 23.389 & 27 & 23.082 \\
  \bottomrule
 \end{tabular}
\caption{Non-restarted FGMRES iteration count for solving the 2D Helmholtz for the Marmousi model.
The RS-CGC 3-grid W(1,1) or V(1,1) cycle with bicubic intergrid operators serves as a preconditioner.
The complex shift used for the combined method is $0.03k^2$.
The time measurements are averaged over 6 experiments.
}
\label{tab:Marmousi}
\end{table}


\subsection{Three-dimensional experiments}\label{subsec:3d}

\begin{figure}
  \begin{center}
    \includegraphics[width=0.6\textwidth, trim=2cm 11cm 2cm 11cm,clip]{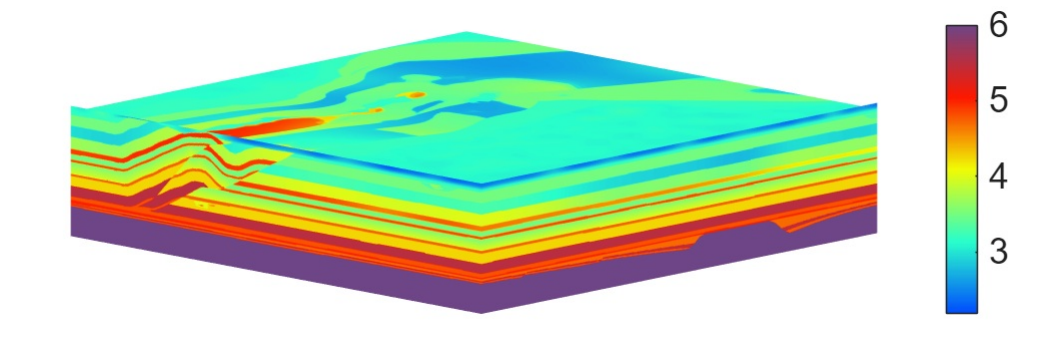}\\
  \caption{The Overthrust model. Velocity units: $km/sec$.}\label{fig:overthrust}
  \end{center}
\end{figure}

The multigrid framework for the 3D experiments is similar to the 2D experiments, except for two differences: 
first, the damping parameters for the Jacobi smoother are $[0.6,0.4]$ for the first and the second grid relaxation respectively.
This values were found in \cite{yovel2025vanka} to yield optimal convergence.
Second, the intergrid operators follow the level-dependent scheme (see Section \ref{sec:main}) for memory considerations.
For the Galerkin coarse operators, we use the real shifts from Table \ref{tab:disp_err_3D}.

\begin{table}
\centering
\begin{tabular}{cc|ccc}
\hline
  \toprule
  \mc{5}{c}{Iteration count, 3D overthrust media}\\
 \midrule
Grid size (cells) & \# DOF & $G=12$ & $G=10$ & time (sec) \\
  \midrule
$64\times64\times32$  & $1.39\cdot10^5$ & 10 & 10 & 4.84 \\
$96\times96\times40$ & $3.86\cdot10^5$ & 10 & 10 & 12.08 \\
$128\times128\times56$  & $9.48\cdot10^5$ & 10 & 10 & 30.00 \\
$160\times160\times64$  & $1.68\cdot10^6$ & 10 & 10 & 42.26 \\
$192\times192\times72$  & $2.72\cdot10^6$ & 10 & 10 & 67.25 \\
  \bottomrule
 \end{tabular}
\caption{Non-restarted GMRES iteration count and wall-clock time for the solution of the 3D Helmholtz equation in overthrust media using RS-CGC preconditioned FGMRES. The time measurements were averaged over 6 tests.}
\label{tab:3d_overthrust}
\end{table}

We apply our method to the Overthrust model \cite{aminzadeh19973}, a geophysical 3D model with nonsmooth coefficients, see Figure \ref{fig:overthrust}. 
The model is defined on a shallow domain of size $20km\times20km\times4.65km$, so we extend the grid downwards by 12--16 points to enable the application of the absorbing boundary layer without adding too much attenuation to the original domain.
We locate the point source in the center of the upper face of the domain, and in this face we do not apply the absorbing layer.

The results in Table \ref{tab:3d_overthrust} show that the iteration count is constant and equals 10, for both $G=12$ and $G=10$, with no additional complex shift.
The average time growth factor is 
$~0.71$, and the time measurements demonstrate linear time scalability.
The accumulation of dispersion error, or the pollution effect, corresponds to the number of cells \emph{per direction}, and the $n_{crit}$ range in Table \ref{tab:disp_err_3D} gives an assesment to the largest number of grid points for which the wavenumber independent convergence that we demonstrate will still hold.

\section{Conclusion}\label{sec:conclusion}

We presented a multigrid solver with a real-shifted coarse grid correction (RS-CGC) which achieves a scalable multigrid method for high frequency Helmholtz problems.
Our solver differs from existing dispersion-correction methods:
the real shift is interwoven into the coarsest Galerkin projection with high grid-transfer operators. 
We presented a grid-to-grid dispersion analysis method that enables to optimize theoretically the real shift of the coarse grid, and to asses the achievable grid sizes for which the method remains scalable.
RS-CGC achieves nearly constant or very low iteration count for homogeneous and heterogeneous media in 2D and 3D, and outperforms the classical CSLP, for 11--12 grid points per wavelength on the fine level.
For 10 grid points per wavelength, we suggested a combination of RS-CGC and CSLP that scales better than each of the methods, and outperforms classical CSLP by an order of magnitude.
We also demonstrated that, despite the large coarse stencils, our method is faster than a re-discretization dispersion-corrected method with standard components.
Despite the inherent limitation underlying any dispersion correction method, the pollution error accumulates slowly enough to enable  large grids, and we provided a rule of thumb to bound these grid sizes.
In future, one can improve this rule of thumb and develop tight bounds for achievable grid sizes.
Another future direction is the application of our method to the elastic Helmholtz equation, through the block-preconditioner we suggested in \cite{yovel2024block}.
Overall, our RS-CGC is a multigrid method with wavenumber independent convergence for high-frequency Helmholtz problems on large heterogeneous media in 2D and in 3D.

\bibliographystyle{siamplain}
\bibliography{RealShiftedSISC.bbl}
\end{document}